\documentclass[12pt]{article}
\usepackage{amssymb}

\newtheorem{thm}{Theorem}[section]
\newtheorem{problem}[thm]{Problem}
\newtheorem{proposition}[thm]{Proposition}
\newtheorem{example}[thm]{Example}

\begin{document}

\title{\vspace{-6ex}~~\\Graph labeling games}

\author{Zsolt Tuza \thanks{
Research supported in part by the National Research,
 Development and Innovation Office -- NKFIH under the grant
  SNN 116095.
  } 
  }
\date{\small Alfr\'ed R\'enyi Institute of Mathematics\\
 \small Hungarian Academy of Sciences\\ \small Budapest, Hungary\\
 \small and\\
  \small Department of Computer Science and Systems Technology\\
    \small University of Pannonia\\ \small Veszpr\'em, Hungary}

\maketitle

\vspace{-3ex}

\begin{abstract}
We propose the study of many new variants of two-person
 graph labeling games.
Hardly anything has been done in this wide open field so far.

~~ 
\end{abstract}

\section{Introduction}\label{intro}

The study of combinatorial games is a classical area in both
 discrete mathematics and game theory; see, e.g., \cite{BCG}.
Nevertheless, until now almost nothing is known about games
 related to graph labeling.
Our goal with this note is to invite attention to this neglected area,
 which offers a wide open field for future research.
Proofs will be published elsewhere.

\subsection{General setting}

We consider labelings of graphs $G=(V,E)$, where $V$ is the vertex set
 and $E$ is the edge set.
The set $L$ of labels will be $\{1, \dots, s\}$ in most cases,
  except in Sections \ref{ss:d1} and \ref{ss:dist}, where we
 shall take $L = \{0, 1, \dots, s\}$, because in the context there,
 the difference between the largest and smallest label in a
 feasible labeling is an important parameter.
We also propose to study further kinds of labels, other than intervals
 of nonnegative integers, but details will not be elaborated here.

Generally speaking, a labeling is a mapping $\phi$ from a certain domain
 into the set $L$ of labels.
In this way we may consider vertex labeling $\phi:V\to L$,
 edge labeling $\phi:E\to L$, total labeling $\phi:(V \cup E)\to L$,
 labeling of faces of a graph on a surface (e.g.\ of a planar graph
 embedded in the plane), etc.

A specific property characteristic for the problem under consideration
 is required for $\phi$; it usually is described in terms of another
 labeling $\phi^*$ derived from $\phi$ by an algebraic operation,
 e.g.\ addition or subtraction.

In the games considered below, two players --- whom we shall call
 Alice and Bob --- alternately select and label vertices or edges
 (typically one vertex or edge in each step) in a graph
 $G=(V,E)$ which is completely known for both players.
The first move is made by Alice, unless stated otherwise.

\subsection{Earlier papers on labeling games}

It seems that the rich literature of graph labeling
 (including more than two thousand works) contains only
 as few as four~(!) papers on labeling games.
Three of them deal with magic labelings \cite{HR,BHSW,GK}, and
 one of them considers the game version of so-called
 $L(d,1)$-labelings \cite{CHKLX}.
  Nothing else seems to have been published in this area so far.
In fact there is a track of research concerning `game chromatic number',
 cf.\ e.g.\ the recent survey \cite{TZ}; but it is considered as
 part of graph coloring, rather than graph labeling.

\subsubsection{Vertex-magic edge- and total labeling games}

 The following game was introduced in 2003 by Hartnell and Rall \cite{HR}.
Given a graph $G=(V,E)$, let the label set be
 $L=\{1, \dots, |E|\}$.
Alice and Bob alternately label a previously unlabeled edge of their
 free choice with a previously unused label.
  A vertex is said to be \emph{full} if all its incident edges are labeled.
The weight of a full vertex is defined as the sum of labels
 of all incident edges.

At the moment when the first full vertex occurs, the magic constant $k$
 is defined as the weight of this vertex.
The basic rule of the game is that all full vertices must have
 the same weight, namely $k$.
In particular, if $vw$ is an unlabeled edge whose selection would make
 both $v$ and $w$ full, but the current edge sums at $v$ and $w$ are
 different, then the players are not allowed to select this edge
 (neither now, nor in any later move).
In general, at any stage of the game, a move is legal if it
 does not create different weights for full vertices.
The winner of the game is the player who makes the last legal move.

Hartnell and Rall \cite{HR} proved that if there are several
 pendant edges incident with each vertex of degree greater than 1,
 then the first player has a winning strategy.
Complementing (in a sense) this result, Giambrone and King
 designed winning strategies on graphs of minimum degree
 at least two.

The `total labeling' version of this game was considered in a single
 paper by Boudreau et al.\ \cite{BHSW}.
Here the label set is $L=\{1, \dots, |V|+|E|\}$, and in each step
 the player on turn labels a previously unlabeled edge or vertex
 with a previously unused label.
The basic rule remains the same: it is not legal to create vertices
 where the sum of all incident edge labels plus the label of the
 vertex itself is not the same.
Similarly to \cite{HR}, the paper \cite{BHSW} presents
 winning strategies for graphs having lots of leaf vertices.

\subsubsection{The $L(d,1)$-labeling game}
  \label{ss:d1}

 The following game was introduced in 2012 by Chia et al.\ \cite{CHKLX}.
Let the label set be $L = \{0, 1, 2, \dots, s\} $,
 the word  ``span'' motivating the notation $s$.
Alice and Bob alternately assign a label  $\phi(v) \in L$ to a
 previously unlabeled vertex $v\in V$.
Labels may be repeated, but the following two rules
 have to be respected:
\begin{itemize}
 \item
  if $vw \in E$, then $|\phi(v) - \phi(w)| \ge d$;
 \item
  if the distance of $v$ and $w$ is 2, then $\phi(v) \ne \phi(w)$.
\end{itemize}

Instead of
 performing the last legal move, Alice's goal
 now
  is to make the entire $G$ labeled, while Bob wants to prevent this.
The problem is to determine the smallest value of the span $s$
 for which Alice has a winning strategy on $G$.
Chia et al.\ \cite{CHKLX} determine an exact formula for the
 cases where $G$ is the complete graph $K_n$ (the minimum
   is $(4d - 1) \cdot n/3 + O(d)$) or the complete
 bipartite graph $K_{p,q}$ plus any
  isolated vertices
 (then the minimum
  is about $p + q + 2d$).

\subsection{Achievement/$\!$Avoidance and Maker--Breaker games}

The two examples above have different characteristics with respect
 to the criteria of winning.
  We next put them in a more general frame separately.
As a general common rule, we shall assume that passing is not allowed.

\subsubsection{Achievement and Avoidance}

The two terms in the subtitle express two opposite ways of how the winner
 of the same game is determined, but on the other hand the two players
 have the same goal.
As an example, in the vertex-magic games defined previously, there is a set of
 winning positions, and the player reaching any of them wins the game.
  That is, the last legal move wins.
In the theory of positional games this is called `normal play';
 adopting the terminology of \cite{HT} we shall call this an
 \emph{Achievement} game since both players aim at achieving
 to move last and reach a final (winning) position.
In the above example the winning positions are the non-extendable
 partial labelings in which all (possibly zero)
  full vertices have the same weight.
Normal play is opposed to `mis\`ere play', what means that the
 last legal move loses.
Again from \cite{HT}, we shall call this an
 \emph{Avoidance} game since both players aim at avoiding
 to move last and reach a final position.
(The analysis of mis\`ere play is extremely difficult in some games.)

In connection with the vertex-magic games introduced above,
 the following natural problems arise.
Concerning \ref{P1-1} and \ref{P1-3} the papers cited above give some
 relevant results, but nothing is known about Problem \ref{P1-2}
  so far.

\begin{problem}   \label{P1-1}
 Describe further classes of graphs on which Alice (or Bob)
 has a winning strategy in the vertex-magic edge- and total
 labeling Achievement games.
\end{problem}

\begin{problem}   \label{P1-2}
 Analyze the Avoidance versions of the vertex-magic edge- and total
 labeling games.
\end{problem}

\begin{problem}   \label{P1-3}
 What is the complexity of determining the winner on\/ $G$ if the
  graph belongs to a specific class?
\end{problem}

\subsubsection{Maker--Breaker}
  \label{ss:dist}

In such games, also including the $L(d,1)$-game,
 the two players have opposite goals.
Alice --- who plays the role of the Maker --- wants to construct
 a specific object
  (in our example it is an $L(d,1)$-labeling of $G$
   with a given value of $s$)
 while Bob --- the Breaker --- wants to prevent this.

\begin{problem}
 Determine the minimum value of\/ $s=s(G)$ for graphs\/ $G$
  from further classes, such that Alice has a winning strategy
  in the\/ $L(d,1)$-labeling on\/ $G$ when the set of labels
  is\/ $\{0,1,\dots,s\}$.
\end{problem}

\begin{problem}
 What is the complexity of determining\/ $s(G)$?
\end{problem}

\begin{problem}
 Study the analogous problems for games on further classes
  of graphs and with other types of
  distance-constrained labelings.
\end{problem}

A general set of constraints for vertex pairs within
 distance $D$ can be described with a $D$-tuple
 $(j_1,j_2,\dots,j_D)$ where $j_1\ge j_2\ge\cdots\ge j_D$.
Then a labeling $\phi:V\to L$ on a graph $G=(V,E)$ is
 feasible if the following property is satisfied:
 $$
   \mathit{if} \ dist(v,w) = i \le D , \
     \mathit{then} \ |\phi(v)-\phi(w)| \ge j_i .
 $$
An interesting particular case is \emph{radio labeling}, in which
 $D$ is taken to be the diameter of $G$, and $j_i = D+1-i$
 for all $1\le i\le diam(G)$.

Apart from $L(d,1)$, no such distance-constrained labeling games
 seem to have been introduced before.


\section{New games from graceful labeling: \hfill\break
 edge-distinguishing games with differences}

In the games proposed in this section, Alice and Bob
 alternately assign a previously unused label
 $\phi(v) \in L=\{1, \dots, s\}$ to a previously unlabeled vertex $v$
 of a given graph $G = (V, E)$.
If both ends of an edge $vw\in E$ are already labeled, then the
 weight of the edge is defined as $|\phi(v) - \phi(w)|$.
  A move is legal if, after it, all edge weights are distinct.
We call such a game \emph{edge-difference distinguishing}.

In the Maker--Breaker version Alice wins if the
 entire $G$ is labeled, and Bob wins if he can prevent this.

\begin{problem}
 Given\/ $G=(V,E)$, for which values of\/ $s$ can Alice win the
  edge-difference distinguishing Maker--Breaker game?
\end{problem}

\begin{problem}
 If Alice can win the edge-difference distinguishing Maker--\break Breaker
  game on\/ $G$ with label set\/ $L=\{1, \dots, s\}$, can she also
  win with\/ $L=\{1, \dots, s+1\}$?
\end{problem}

\begin{problem}
 What is the complexity of deciding whether Alice can win
  on an input graph\/ $G$ with label set\/ $L=\{1, \dots, s\}$?
\end{problem}

\begin{problem}
 Study the Achievement and the Avoidance versions of this game.
\end{problem}

If $|V|$ is odd and Alice has a winning strategy in the Maker--Breaker
 version, then of course she can also win the Achievement version.
But other implications do not seem to be obvious.

\subsection{The graceful game}

We next present very simple examples for the case where
 $s=|E|+1$ and hence Alice's goal is to end up with a
 graceful labeling of $G$.
Instead of `edge-difference distinguishing' we may simply
 call it the \emph{Graceful Game}.


\begin{proposition}
 Alice can win the Graceful Game on any star\/  $K_{1,n-1}$;
   and Bob has a winning strategy on
  every path of order at least four.
\end{proposition}

\begin{proposition}
 Alice can win the Graceful Game on complete graphs with at most
  three vertices, no matter which player starts;
  and Bob can win on every larger complete graph.
\end{proposition}

\section{New games from antimagic labeling: \hfill\break
 edge-distinguishing games with sums}

In the \emph{edge-sum distinguishing games} proposed in this section,
 Alice and Bob alternately assign a previously unused label
 $\phi(v) \in L=\{1, \dots, s\}$ to a previously unlabeled vertex $v$
 of a given graph $G = (V, E)$.
If both ends of an edge $vw\in E$ are already labeled, then the
 weight of the edge is defined as $\phi(v) + \phi(w)$.
  A move is legal if, after it, all edge weights are distinct.

In the Maker--Breaker version Alice wins if the
 entire $G$ is labeled, and Bob wins if he can prevent this.

\begin{problem}
 Given\/ $G=(V,E)$, for which values of\/ $s$ can Alice win the
  edge-sum distinguishing Maker--Breaker game?
\end{problem}

\begin{problem}
 If Alice can win the edge-sum distinguishing Maker--\break Breaker
  game on\/ $G$ with label set\/ $L=\{1, \dots, s\}$, can she also
  win with\/ $L=\{1, \dots, s+1\}$?
\end{problem}

\begin{problem}
 What is the complexity of deciding whether Alice can win
  on an input graph\/ $G$ with label set\/ $L=\{1, \dots, s\}$?
\end{problem}

\begin{problem}
 Study the Achievement and the Avoidance versions of this game.
\end{problem}

\subsection{The game on cycles}

 In case of $G=C_n$ we have $|V|=|E|=n$, and now take $s=n$.
On $C_3$ the game is trivial win for Alice because eventually
 a cyclic permutation of
 $\{1,2,3\}$ will occur, no matter how the players play.

\begin{proposition}
 Bob can win the game on each of\/ $C_4$,\/ $C_5$, and\/ $C_6$.
\end{proposition}

Modifying the rules by disregarding the distinct sums condition,
 for every even $n\ge 4$, in the Bob-start game Alice can achieve a
 complete final labeling in which the sum $n+1$ does not occur at all.
Indeed, if Bob assigns label $l_i$ to a vertex in his $i$th move,
 Alice can simply assign $n+1-i$ to the corresponding antipodal vertex.

\section{Variants of games}

There are many details where one can make his/her favorite
 choice, each combination defining a different game.
Instead of proposing further games explicitly, we only list
 here some aspects to be taken into consideration.

\begin{itemize}
 \item
  Label the vertices, or the edges, or both, or other parts
   of $G$ (e.g., subgraphs isomorphic to a specified graph $F$,
   or the faces if $G$ is embedded in the plane or
   on another surface).
 \item
  Achievement, or Avoidance; or Maker--Breaker.
 \item
  Biased game, frequently called $(a:b)$-game in the literature
   --- for specified $a,b\in \mathbb{N}$, Alice makes
   $a$ consecutive moves, followed by Bob's $b$ consecutive moves,
   and they alternate in this way until the game terminates.
 \item
  Labels may or may not be repeated.
 \item
  Types of labels used in the labeling
   (integers, natural/rational/real/complex numbers,
    algebraic variables, vectors, subsets, multisets, ...)
 \item
  Type of labeling --- see \cite{G}.
 \item
  Way of computing $\phi^*$ from $\phi$ (sum, difference, product, mod $q$,
   vector- and (multi)set-operations, ...)
 \item
  Local conditions or global ones.
 \item
  `Legal move' means the same for both players, or some of the
   restrictive rules have to be respected only by Alice
   (or only by Bob).
 \item
  Passing is forbidden, or one or both players are allowed to pass.
 \item
  Which of the players starts the game.
\end{itemize}

We should note that in some variants
 additional restrictive rules may be reasonable.
Consider, for instance, the following.

\begin{example}
 In an edge-distinguishing Maker--Breaker game the vertices of
  $G=(V,E)$ get labels from $L=\{1, \dots, s\}$, but labels may be
  repeated; all `full' edges $vw$ (with both ends labeled) get
   weight $f(\phi(v),\phi(w))$, where $f$ is a specified symmetric
   two-variable function (e.g., $f(x,y)$ is $x+y$ or $|x-y|$).
 A move is legal if, after it, any two \emph{adjacent} full edges
  have different weights.
\end{example}

If $G$ is connected and does not have a
 vertex
  adjacent to all
 the others,
   then
 Bob
 has
  a trivial winning
 way:
After Alice's first move, he
 assigns
 the same label to
 a vertex at distance two.
Then the two edges to the common neighbor
 will
 have the same weight.
To exclude this trivial
 winning
 it is reasonable to require
 that vertices with identical labels be
 at distance at least 3 apart.

\section{Conclusion}

Until now only very few results are known on game versions of
 graph labeling.
In this note we propose to study further ones from the very
 rich collection of variations.
It will be a subject of future research to identify those
 games which lead to really interesting results.
We should mention, however, that the games discussed above
 motivated some coloring games, too; cf.\ \cite{T}.

Beyond graphs, many more combinatorial structures may also
 turn out to be interesting in this direction as well ---
 directed graphs, multigraphs, hypergraphs, tournaments,
 partially ordered sets, ...

\end{document}